\documentclass[12pt]{article}
\usepackage{amsmath,amssymb,amsbsy,amsfonts,amsthm,latexsym,
            amsopn,amstext,amsxtra,euscript,amscd}
\begin{document}
\bibliographystyle{plain}
\newfont{\teneufm}{eufm10}
\newfont{\seveneufm}{eufm7}
\newfont{\fiveeufm}{eufm5}
%
%
\newfam\eufmfam
              \textfont\eufmfam=\teneufm \scriptfont\eufmfam=\seveneufm
              \scriptscriptfont\eufmfam=\fiveeufm
\def\bbbr{{\rm I\!R}}
\def\bbbm{{\rm I\!M}}
\def\bbbn{{\rm I\!N}}
\def\bbbf{{\rm I\!F}}
\def\bbbh{{\rm I\!H}}
\def\bbbk{{\rm I\!K}}
\def\bbbp{{\rm I\!P}}
\def\bbbone{{\mathchoice {\rm 1\mskip-4mu l} {\rm 1\mskip-4mu l}
{\rm 1\mskip-4.5mu l} {\rm 1\mskip-5mu l}}}
\def\bbbc{{\mathchoice {\setbox0=\hbox{$\displaystyle\rm C$}\hbox{\hbox
to0pt{\kern0.4\wd0\vrule height0.9\ht0\hss}\box0}}
{\setbox0=\hbox{$\textstyle\rm C$}\hbox{\hbox
to0pt{\kern0.4\wd0\vrule height0.9\ht0\hss}\box0}}
{\setbox0=\hbox{$\scriptstyle\rm C$}\hbox{\hbox
to0pt{\kern0.4\wd0\vrule height0.9\ht0\hss}\box0}}
{\setbox0=\hbox{$\scriptscriptstyle\rm C$}\hbox{\hbox
to0pt{\kern0.4\wd0\vrule height0.9\ht0\hss}\box0}}}}
\def\bbbq{{\mathchoice {\setbox0=\hbox{$\displaystyle\rm
Q$}\hbox{\raise
0.15\ht0\hbox to0pt{\kern0.4\wd0\vrule height0.8\ht0\hss}\box0}}
{\setbox0=\hbox{$\textstyle\rm Q$}\hbox{\raise
0.15\ht0\hbox to0pt{\kern0.4\wd0\vrule height0.8\ht0\hss}\box0}}
{\setbox0=\hbox{$\scriptstyle\rm Q$}\hbox{\raise
0.15\ht0\hbox to0pt{\kern0.4\wd0\vrule height0.7\ht0\hss}\box0}}
{\setbox0=\hbox{$\scriptscriptstyle\rm Q$}\hbox{\raise
0.15\ht0\hbox to0pt{\kern0.4\wd0\vrule height0.7\ht0\hss}\box0}}}}
\def\bbbt{{\mathchoice {\setbox0=\hbox{$\displaystyle\rm
T$}\hbox{\hbox to0pt{\kern0.3\wd0\vrule height0.9\ht0\hss}\box0}}
{\setbox0=\hbox{$\textstyle\rm T$}\hbox{\hbox
to0pt{\kern0.3\wd0\vrule height0.9\ht0\hss}\box0}}
{\setbox0=\hbox{$\scriptstyle\rm T$}\hbox{\hbox
to0pt{\kern0.3\wd0\vrule height0.9\ht0\hss}\box0}}
{\setbox0=\hbox{$\scriptscriptstyle\rm T$}\hbox{\hbox
to0pt{\kern0.3\wd0\vrule height0.9\ht0\hss}\box0}}}}
\def\bbbs{{\mathchoice
{\setbox0=\hbox{$\displaystyle     \rm S$}\hbox{\raise0.5\ht0\hbox
to0pt{\kern0.35\wd0\vrule height0.45\ht0\hss}\hbox
to0pt{\kern0.55\wd0\vrule height0.5\ht0\hss}\box0}}
{\setbox0=\hbox{$\textstyle        \rm S$}\hbox{\raise0.5\ht0\hbox
to0pt{\kern0.35\wd0\vrule height0.45\ht0\hss}\hbox
to0pt{\kern0.55\wd0\vrule height0.5\ht0\hss}\box0}}
{\setbox0=\hbox{$\scriptstyle      \rm S$}\hbox{\raise0.5\ht0\hbox
to0pt{\kern0.35\wd0\vrule height0.45\ht0\hss}\raise0.05\ht0\hbox
to0pt{\kern0.5\wd0\vrule height0.45\ht0\hss}\box0}}
{\setbox0=\hbox{$\scriptscriptstyle\rm S$}\hbox{\raise0.5\ht0\hbox
to0pt{\kern0.4\wd0\vrule height0.45\ht0\hss}\raise0.05\ht0\hbox
to0pt{\kern0.55\wd0\vrule height0.45\ht0\hss}\box0}}}}
\def\bbbz{{\mathchoice {\hbox{$\sf\textstyle Z\kern-0.4em Z$}}
{\hbox{$\sf\textstyle Z\kern-0.4em Z$}}
{\hbox{$\sf\scriptstyle Z\kern-0.3em Z$}}
{\hbox{$\sf\scriptscriptstyle Z\kern-0.2em Z$}}}}
\def\ts{\thinspace}

\newtheorem{theorem}{Theorem}
\newtheorem{lemma}[theorem]{Lemma}
\newtheorem{claim}[theorem]{Claim}
\newtheorem{cor}[theorem]{Corollary}
\newtheorem{prop}[theorem]{Proposition}
\newtheorem{definition}[theorem]{Definition}
\newtheorem{remark}[theorem]{Remark}
\newtheorem{question}[theorem]{Open Question}

\def\qed{\ifmmode
\squareforqed\else{\unskip\nobreak\hfil
\penalty50\hskip1em\null\nobreak\hfil\squareforqed
\parfillskip=0pt\finalhyphendemerits=0\endgraf}\fi}

\def\squareforqed{\hbox{\rlap{$\sqcap$}$\sqcup$}}

\def \C {{\mathbb C}}
\def \F {{\mathbb F}}
\def \L {{\mathbb L}}
\def \K {{\mathbb K}}
\def \Q {{\mathbb Q}}
\def \Z {{\mathbb Z}}
\def\cA{{\mathcal A}}
\def\cB{{\mathcal B}}
\def\cC{{\mathcal C}}
\def\cD{{\mathcal D}}
\def\cE{{\mathcal E}}
\def\cF{{\mathcal F}}
\def\cG{{\mathcal G}}
\def\cH{{\mathcal H}}
\def\cI{{\mathcal I}}
\def\cJ{{\mathcal J}}
\def\cK{{\mathcal K}}
\def\cL{{\mathcal L}}
\def\cM{{\mathcal M}}
\def\cN{{\mathcal N}}
\def\cO{{\mathcal O}}
\def\cP{{\mathcal P}}
\def\cQ{{\mathcal Q}}
\def\cR{{\mathcal R}}
\def\cS{{\mathcal S}}
\def\cT{{\mathcal T}}
\def\cU{{\mathcal U}}
\def\cV{{\mathcal V}}
\def\cW{{\mathcal W}}
\def\cX{{\mathcal X}}
\def\cY{{\mathcal Y}}
\def\cZ{{\mathcal Z}}
\newcommand{\rmod}[1]{\: \mbox{mod}\: #1}

\def\tcN{\cN^\mathbf{c}}
\def\F{\mathbb F}
\def\Tr{\operatorname{Tr}}
\def\mand{\qquad \mbox{and} \qquad}
\renewcommand{\vec}[1]{\mathbf{#1}}
\def\eqref#1{(\ref{#1})}
\newcommand{\ignore}[1]{}
\hyphenation{re-pub-lished}
\parskip 1.5 mm
\def\lln{{\mathrm Lnln}}
\def\Res{\mathrm{Res}\,}
\def\F{{\bbbf}}
\def\Fp{\F_p}
\def\fp{\Fp^*}
\def\Fq{\F_q}
\def\ff{\F_2}
\def\ffn{\F_{2^n}}
\def\K{{\bbbk}}
\def \Z{{\bbbz}}
\def \N{{\bbbn}}
\def\Q{{\bbbq}}
\def \R{{\bbbr}}
\def \P{{\bbbp}}
\def\Zm{\Z_m}
\def \Um{{\mathcal U}_m}
\def \Bf{\frak B}
\def\Km{\cK_\mu}
\def\va {{\mathbf a}}
\def \vb {{\mathbf b}}
\def \vc {{\mathbf c}}
\def\vx{{\mathbf x}}
\def \vr {{\mathbf r}}
\def \vv {{\mathbf v}}
\def\vu{{\mathbf u}}
\def \vw{{\mathbf w}}
\def \vz {{\mathbfz}}
\def\\{\cr}
\def\({\left(}
\def\){\right)}
\def\fl#1{\left\lfloor#1\right\rfloor}
\def\rf#1{\left\lceil#1\right\rceil}
\def\flq#1{{\left\lfloor#1\right\rfloor}_q}
\def\flp#1{{\left\lfloor#1\right\rfloor}_p}
\def\flm#1{{\left\lfloor#1\right\rfloor}_m}
\def\Al{{\sl Alice}}
\def\Bob{{\sl Bob}}
\def\Or{{\mathcal O}}
\def\inv#1{\mbox{\rm{inv}}\,#1}
\def\invM#1{\mbox{\rm{inv}}_M\,#1}
\def\invp#1{\mbox{\rm{inv}}_p\,#1}
\def\Ln#1{\mbox{\rm{Ln}}\,#1}
\def \nd {\,|\hspace{-1.2mm}/\,}
\def\ord{\mu}
\def\E{\mathbf{E}}
\def\Cl{{\mathrm {Cl}}}
\def\epp{\mbox{\bf{e}}_{p-1}}
\def\ep{\mbox{\bf{e}}_p}
\def\eq{\mbox{\bf{e}}_q}
\def\bm{\bf{m}}
\newcommand{\floor}[1]{\lfloor {#1} \rfloor}
\newcommand{\comm}[1]{\marginpar{
\vskip-\baselineskip 
\raggedright\footnotesize
\itshape\hrule\smallskip#1\par\smallskip\hrule}}
\def\rem{{\mathrm{\,rem\,}}}
\def\dist {{\mathrm{\,dist\,}}}
\def\etal{{\it et al.}}
\def\ie{{\it i.e. }}
\def\veps{{\varepsilon}}
\def\eps{{\eta}}
\def\ind#1{{\mathrm {ind}}\,#1}
               \def \MSB{{\mathrm{MSB}}}
\newcommand{\abs}[1]{\left| #1 \right|}

\title {On the Counting Function of Elliptic Carmichael Numbers}
%

\author{{\sc Florian~Luca}\\ 
Centro de Ciencias Matem{\'a}ticas,\\
Universidad Nacional Autonoma de M{\'e}xico,\\
C.P. 58089, Morelia, Michoac{\'a}n, M{\'e}xico\\
{\tt fluca@matmor.unam.mx}
\and
{\sc Igor E. Shparlinski}  \\
Deptartment of Computing  \\
Macquarie University \\
Sydney, NSW 2109, Australia\\
{\tt igor@ics.mq.edu.au}
}

\date{\today}
\pagenumbering{arabic}

\maketitle

\begin{abstract} We give an upper bound for the number elliptic Carmichael numbers $n \le x$ 
that have recently been introduced by J.~H.~Silverman. We also discuss
several possible ways for further improvements. 
\end{abstract}

\section{Introduction}

Let $E$ be an elliptic curve over the field of rational numbers $\Q$ given by 
an {\it affine Weierstra\ss\ equation}:
$$
E~:~Y^2 = X^3 + aX + b.
$$
In particular, it has a nonzero discriminant  $\Delta = 4a^3+27b^2$. 
We  refer to~\cite{Silv1} for a background on elliptic curves. 

For a prime $p$,  we define $a_p$ by $\# E(\F_p) = p+1 - a_p$,
where $E(\F_p)$ in the set of $\F_p$-rational points on the
reduction of $E$ modulo $p$ including the point at infinity $O_p$.  

We also recall that if $p\nmid \Delta$, then $E(\F_p)$ has a
structure of an Abelian group (see~\cite[Chapter~III, Section~2]{Silv1}). 

Since by  the {\it Hasse bound\/} $a_p = O(p^{1/2})$ 
(see, for example,~\cite[Chapter~V, Theorem~1.1]{Silv1}),
for $\Re s > 3/2$ we can define the  $L$-function
$$
L(s) = \prod_{p\nmid \Delta} \(1-a_p p^{-s}\)^{-1}
\prod_{p\mid \Delta} \(1-a_p p^{-s}+p^{1-2s}\)^{-1},
$$
which we expand to the power series
$$
L(s)  = \sum_{n=1}^\infty \frac{a_n}{n^s}
$$
(see, for example,~\cite[Chapter~V, Exercise~8.19]{Silv1}).

Slightly relaxing the definition given in~\cite{Silv2} and thus expanding the class of 
numbers we consider, we say that a positive integer $n$ is an {\it $E$-Carmichael number\/}
if 
\begin{itemize}
\item it is not a prime power;
\item for any prime divisor  $p\mid n$ we have $p\nmid \Delta$;
\item for any point $P \in E(\F_p)$ we have
\begin{equation}
\label{eq:E-Carm}
(n+1-a_n) P = O_p,
\end{equation}
where both the equation and the group law are considered over $\F_p$.
\end{itemize}

Here we show that the sequence $E$-Carmichael numbers is 
of asymptotic density zero.

\section{Notation}

We recall that the notations $U = O(V)$, $U \ll V$ and  $V \gg U$  are all
equivalent to the statement that the inequality $|U| \le c\,V$ holds
with some constant $c> 0$. 
Throughout the paper, any implied constants
in the symbols `$O$',  `$\ll$' and $\gg$' may occasionally depend, where obvious, on
the curve $E$,  and are absolute otherwise.

 We write $\log_1 x=\max\{1,\log x\}$. For an integer $k\ge 2$, we write $\log_k x$ for the iteratively defined function 
given by $\log_k x=\log_1 (\log_{k-1} x)$. When $k=1$ we omit the subscript and thus understand that all natural logarithms that appear exceed $1$.

\section{Main Result}

For a real $x \ge 1$,  let $N_E(x)$ be the number of $E$-Carmichael numbers $n \le x$.

\begin{theorem}
\label{thm:ECCarm}  For a sufficiently large $x$
$$
N_E(x) \ll x\frac{(\log_3 x)^{1/2} (\log_4 x)^{1/2}}{(\log_2 x)^{1/4}}.
$$
\end{theorem}

\section{Preparations}

We start with an integer $a \ne 0, \pm 2$ and a special case of a result Serre~\cite{Ser} that gives an upper bound on 
$$
\pi_{E}(x;a)=\#\{p\le x~:~a_p = a\}.
$$

\begin{lemma}
\label{lem:Ser}
The estimate 
$$
\pi_E(x;a)\ll \pi(x)\frac{ (\log_2 x)^{2/3}(\log_3 x)^{1/3}}
{(\log x)^{1/3}}
$$
holds for all $a\ne 0,\pm 2$, where the implied constants depend only on the elliptic curve $E$.
\end{lemma}

We also need the following result of David and Wu~\cite[Theorem~2.3~(i)]{DaWu}, which improves and generalises several 
previous bounds (see~\cite{CoFoMu,CoLuSh}). For integers $a$ and $b\ge 1$ let 
$$
\pi_{E}(x;a,b)=\#\{p\le x~:~\#E(\F_p)\equiv a\pmod b\}.
$$

Let $\varphi(k)$ denote  the Euler function of an integer $k\ge 1$. 

\begin{lemma}
\label{lem:DaWu}
The estimate 
$$
\pi_E(x;a,b)\ll \frac{\pi(x)}{\varphi(b)}+x\exp\left(-A b^{-2} {\sqrt{\log x}}\right)
$$
holds uniformly for $\log x\gg b^{12}\log b$, where the implied constants depend only on the elliptic curve $E$ and $A$ is a positive absolute constant.
\end{lemma}

\section{Proof of Theorem~\ref{thm:ECCarm}}.

Let $t_p$ be the exponent of the group $E(\F_p)$,
that is, the largest possible order  of any point $P \in E(\F_p)$. 

We see from~\eqref{eq:E-Carm} that for any $E$-Carmichael number $n$
we have 
\begin{equation}
\label{eq:tpn}
t_p \mid n+1-a_n
\end{equation}
for all primes $p\mid n$.

Now fix some $z > y > 1$ and remove $n\le x$ without a
prime divisor in $[y,z]$. Let $\cE_1(x)$ be the set of such $n$. 
By the Brun sieve, see~\cite[Section~I.4.2]{Ten} and Mertens' formula, 
see~\cite[Section~I.1.6]{Ten}, we have
\begin{equation}
\label{eq:E1}
\#\cE_1(x)\ll x\prod_{y\le p\le z} \left(1-\frac{1}{p}\right)=O\left(x\frac{\log y}{\log z}\right).
\end{equation}
Then remove all $n\le x$ such that $p^2 \mid n$ for some
$p\ge y$. Let $\cE_2(x)$ be the set of such $n$. Fixing $p$, the number of $n\le x$ which are divisible by $p^2$ is at most $x/p^2$. Hence,
\begin{equation}
\label{eq:E2}
\#\cE_2(x)\le \sum_{y\le p\le z} \frac{x}{p^2}=O\left(\frac{x}{y}\right).
\end{equation}
Let $P(n)$ be the largest prime factor of $n$. We remove $n\le x$ such that $P(n)\le w$, 
where
$$
w=\exp\left(\frac{\log x \log_4 x}{2\log_3 x}\right).
$$
Put $\cE_3(x)$ for the set of such $n$. It is well-known that 
$$
\#\cE_3(x)=\frac{x}{\exp((1+o(1))u\log u)},
$$
as $x\to\infty$, where
$$
u=\frac{\log x}{\log w}=\frac{2\log_3 x}{\log_4 x}.
$$
Since
$$
u\log u=(2+o(1))\log_3 x
$$
as $x\to\infty$, we derive
\begin{equation}
\label{eq:E3}
\#\cE_3(x)=\frac{x}{(\log_2 x)^{2+o(1)}}=O\left(\frac{x}{\log_2 x}\right).
\end{equation}

Assume that $w^{1/2}>2z$. Then any remaining integer $n\le x$ can be written under the form $n = pPm$, where $p\in [y,z],~P=P(n)>w$ and 
$pP$ is coprime to $m$. Since the coefficient $a_n$ is a multiplicative function of $n$, we have  $a_n = a_m a_pa_P.$
Then, we see from~\eqref{eq:tpn}, that
\begin{equation}
\label{eq:tpmp}
t_p | mPp + 1 - a_m a_pa_P.   
\end{equation}
We fix $p\in [y,z]$  count the number of choices for the pair $(m,P)$. Assume next that $p\mid t_p$. Let $\cE_4(x)$ be the number of such $n$. In this case, $t_p=p$, $a_p=1$ and 
congruence~\eqref{eq:tpmp} shows that $p\mid a_{mP}$. 

Estimating the number of such products $mP\le x/p$ trivially  as $O(x/p)$, 
summing up over all $p\in [y,z]$ with $a_p=1$ and 
using Abel's summation formula and Lemma~\ref{lem:Ser}, we derive 
\begin{equation}
\label{eq:E4}
\# \cE_4(x)\ll \sum_{\substack{y\le p\le z\\ a_p = 1}}\frac{x}{p} 
\ll \frac{ x(\log_2 y)^{2/3}(\log_3 y )^{1/3}}{(\log y)^{1/3}}.
\end{equation}

From now on, we assume that $t_p$ and $p$ are coprime. 
Note that $t_p \gg p^{1/2}$ (see~\cite{Schoof} for a  slightly more 
precise result). We next write
$$
t_p=d_1d_2,
$$
where $d_1=\gcd(t_p,m)$. Suppose that $d_1>t_p^{1/2}$ and let $\cE_5(x)$ be the set of such $n\le x$. Then $m=d_1m_1$, so $n$ is a multiple of $pd_1$. The number of such choices when $p$ and $d_1\mid t_p$ are  fixed is at most $x/pd_1=O(x/p^{5/4})$. Summing up over all primes $p$ and divisors $d_1$ of $t_p$ which exceed $t_p^{2/3}$, we get that
\begin{equation}
\label{eq:E5}
\#\cE_5(x)\ll \sum_{y\le p\le z} \frac{\tau(t_p)}{p^{5/4}}=O\left(\frac{x}{y^{1/4+o(1)}}\right)
\end{equation}
as $y \to \infty$. 

Let $\cE_6(x)$ be the set of the remaining $n\le x$. Writing again $m=d_1 m_1$, the divisibility relation~\eqref{eq:tpmp} implies that $d_1\mid a_p a_m a_P$. Fix also $m$ and we put $d_3=\gcd(d_1,a_p),~d_4=\gcd(d_1/d_3,a_m)$, and $d_5=d_1/(d_3d_4)$. Then the relation $a_P=d_5\lambda$ holds with some positive integer 
$\lambda$. Further, the divisibility relation~\eqref{eq:tpmp} gives
$$
d_2\mid m_1p P-\left(\frac{a_p}{d_3}\right)\left(\frac{a_m}{d_4}\right)\lambda, 
$$
and $m_1 p$ is invertible modulo $d_2$. This shows that
\begin{equation}
\label{eq:congforP}
P\equiv \left(\frac{a_p}{d_3}\right)\left(\frac{a_m}{d_5}\right)\lambda \pmod {d_2}.
\end{equation}
In the right--hand side of the congruence~\eqref{eq:congforP}, we assume that $a_p/d_3$ and $a_m/d_5$ are coprime to $d_2$, otherwise $P\mid d_2$, which is impossible since it would lead to 
$$
w\le P\le d_2\le t_p<p+2{\sqrt{p}}+1<2z,
$$ 
for large $x$, which is impossible. Observe that the value of $\lambda\pmod {d_2}$ determines both $P$ and 
$a_P$ modulo $d_2$. 
In turn, these define $\#E(\F_P)$ modulo $P$. 
By Lemma~\ref{lem:DaWu}, we derive that number of such $P\le x/(mp)$ is of order at most
\begin{equation}
\label{eq:Cheb1}
\begin{split}
&\frac{\pi(x/mp)}{\varphi(d_2)} + \frac{x}{mp}\exp\(-A d_2^{-2} {\sqrt{\log x}}\)\\
& \qquad\qquad \ll
\frac{x}{mp\varphi(d_2)\log(x/mp)} + x\exp\(-Ad_2^{-2} {\sqrt{\log x}}\), 
\end{split}
\end{equation}
provided that 
$$
d_2\log d_2 \le (\log(x/mp))^{1/12}.
$$ 
Since $d_2\le t_p\le 2z$ and $x/mp\ge P\ge w$, so 
$$
\log(x/mp)\ge \log w\ge \frac{ \log x \log_3 x}{\log_2 x} ,
$$ 
it follows that the above inequality holds if we choose 
\begin{equation}
\label{eq:xyz1}
z\le (\log x)^{1/13}
\end{equation}
 and $x$ is sufficiently large.
For such values of $x$ and $z$, the second term in the estimate~\eqref{eq:Cheb1} is 
$$
x\exp\(-Ad_2^{-2} {\sqrt{\log x}}\) \le x \exp\(-0.25A (\log x)^{11/26}\),
$$
and is negligible compared with the first. So, the number of such primes $P\le x/(mp)$ is of order at most
$$
 \frac{x}{mp\varphi(d_2) \log(x/mp)}\ll \frac{x\log_2 z}{mp d_2 \log(x/mp)},
$$
where we have used that,  by the well-known bound on the minimal order of the Euler function (see~\cite[Section~I.5.4]{Ten}), the lower bound
$$
\varphi(d_2)\gg d_2/\log_2 d_2\gg d_2/\log_2 z
$$
holds for all $d_2\le t_p\le 2z$. 
Since $x/(mp)>w/z>w^{1/2}$ and also since $d_2=t_p/d_1\ge t_p^{1/2}\gg p^{1/4}$, we get that the above estimate is of order at most
$$
\frac{x\log_3 x\log_2 z}{mp d_2 \log x\log_4 x}\ll  \frac{x \log_3 x\log_2 z}{mp^{5/4} \log x\log_4 x}.
$$
Now we sum up the above inequality over all $p\in [y,z]$, all quadruple of divisors $(d_1,d_2,d_3,d_4)$ of $t_p$ and over all $m$ getting a bound of shape
$$
\frac{x\log_3 x\log_2 z }{\log x\log_4 x} \sum_{y\le p\le z} \sum_{m\le x}  \frac{\tau(t_p)^4}{m p^{5/4}}\ll \frac{x\log_3 x\log_2 z}{y^{1/4+o(1)}\log_4 x},
$$
as $x\to\infty$.
Thus, we get that
\begin{equation}
\label{eq:E6}
\#\cE_6(x)\le \frac{x\log_3 x \log_2 z}{y^{1/4+o(1)}\log_4 x}
\end{equation}
as $x\to\infty$. 
From the estimates~\eqref{eq:E1}, \eqref{eq:E2},  \eqref{eq:E3}, \eqref{eq:E4}, \eqref{eq:E5} and~\eqref{eq:E6}, we conclude that
\begin{equation*}
\begin{split}
\#N_E(x)  \ll  x&\left(\frac{\log y}{\log z}+\frac{1}{y}+\frac{1}{\log_2 x}+\frac{(\log_2 y)^{2/3} (\log_3 y)^{1/3}}{(\log y)^{1/3}}+\frac{1}{y^{1/4+o(1)}}
\right.\\
& \qquad\qquad\qquad\qquad\qquad\qquad \qquad\qquad\quad+  \left.
\frac{\log_3 x \log_2 z}{y^{1/4+o(1)} \log_4 x}\right).
\end{split}
\end{equation*}
Since $z\le (\log x)^{1/13}$, the third term is dominated by the first and the second term is dominated by the fourth.  
Since $y\le z\le (\log x)^{1/13}$, it follows that
$$
(\log_2 y)^{2/3} (\log_3 y)^{2/3}\ll (\log_3 x)^{2/3} (\log_4 x)^{2/3},
$$
so we see that
$$
\frac{(\log_2 y)^{2/3} (\log_3 y)^{1/3}}{(\log y)^{1/3}}+\frac{\log_3 x \log_2 z}{y^{1/4+o(1)}\log_4 x}\ll \frac{(\log_3 x)^{2/3} (\log_4 x)^{2/3}}{(\log y)^{1/3}},
$$
provided that 
\begin{equation}
\label{eq:xyz2}
y^{1/5}>(\log_3 x)^2\ge \log_3 x\log_2 z.
\end{equation}
It now follows easily that 
$$
\cN_E(x) \ll x \left(\frac{\log y}{\log z}+\frac{(\log_3 x)^{2/3}(\log_4 x)^{2/3}}{(\log y)^{1/3}}\right).
$$
We now choose 
$$z=(\log x)^{1/14}\mand  
y=\exp\left((1/14)(\log_2 x)^{3/4} (\log_3 x)^{1/2} (\log_4 x)^{1/2}\right),
$$
thus~\eqref{eq:xyz1} and~\eqref{eq:xyz2} are satisfied,  
and we derive the desired result.

\section{Comments} 

We recall that under the Generalised Riemann Hypothesis, Serre~\cite{Ser}
gives a much stronger estimate
$$
\pi_E(x;a)\ll \pi(x) x^{-1/6}(\log  x)^{2/3}, \qquad a\ne 0,\pm 2,
$$
instead of that of Lemma~\ref{lem:Ser}; we also refer  to~\cite{Coj}
for a survey of other results and conjectures related to Lemma~\ref{lem:Ser}. 
Furthermore, also under the Generalised Riemann Hypothesis, 
David and Wu~\cite[Theorem~2.3~(iii)]{DaWu} show that one has 
the estimate 
$$
\pi_E(x;a,b)\ll \frac{\pi(x)}{\varphi(b)}
$$
uniformly for $b \ll  x^{1/8}/\log x$, instead of that of Lemma~\ref{lem:DaWu}.
Using these bounds in our argument, one can easily obtain a conditional 
improvement of  Theorem~\ref{thm:ECCarm}. 
It is also possible that for CM curves one can also obtain stronger 
results. For example, in~\cite{Coj} one can find a survey of improvements of 
Lemma~\ref{lem:Ser} for CM curves. There is little doubt that 
Lemma~\ref{lem:DaWu} can also be improved for CM curves. 
However, in order to get substantially 
better bounds, our argument, which treats 
the elements the set $\#\cE_1(x)$ trivially and relies
on the bound~\eqref{eq:E1}, 
ought to be augmented with some new ideas. 

Another approach to a possible improvement  of  Theorem~\ref{thm:ECCarm}
is via a more efficient treatment of elements of the set $\cE_4(x)$. 
In turn, this leads to a question of obtaining nontrivial upper 
bounds on the cardinality of the set 
$$
\{n\le x~:~ a_n \equiv a \pmod p\}
$$
for a prime $p$ and an integer $a$ (only the case $a=1$ is relevant 
to our applications). Obtaining such bounds is certainly of independent 
interest.

\section{Acknowledgements}

The authors would like to thank Alina Cojocaru for useful discussions. 

During the preparation of this paper,
and F.~L. was supported in part by Project 
PAPIIT IN104512 and a Marcos Moshinsky fellowship and I.~S. by ARC Grant DP1092835 (Australia)
and by NRF Grant~CRP2-2007-03 (Singapore).

\end{document}